\documentclass[11pt]{article}
\usepackage{latexsym}
\usepackage{epsfig}
\usepackage{amsfonts}
\usepackage{amssymb}
\usepackage{amsmath}
\usepackage{amsthm}

\usepackage{color}
\definecolor{purple}{rgb}{0.65, 0, 1}
\definecolor{orange}{rgb}{1,.5,0}
\definecolor{brown}{rgb}{.9,.73,.26}

\def\e{\epsilon}

\def\II{{\rm I\kern-0.5exI}}
\def\I{{\rm I\kern-0.5exI\kern-0.5exI}}

\newtheorem{theorem}{Theorem}[section]
\newtheorem{lemma}[theorem]{Lemma}

\newtheorem{proposition}[theorem]{Proposition}
\newtheorem{corollary}[theorem]{Corollary}

\newcommand{\lp}{\left(}
\newcommand{\rp}{\right)}

\newcommand{\abs}[1]{\left\lvert #1\right\rvert}
\newcommand{\norm}[1]{\left\lVert #1\right\rVert}
\newcommand{\pd}[2]{\frac{\partial #1}{\partial #2}}

\newcommand{\be}{\begin{equation}}
\newcommand{\ee}{\end{equation}}
\newcommand{\bee}{\begin{equation*}}
\newcommand{\eee}{\end{equation*}}
\newcommand{\bea}{\begin{eqnarray}}
\newcommand{\eea}{\end{eqnarray}}
\newcommand{\bs}{\begin{split}}
\newcommand{\es}{\end{split}}

\begin{document}
\title{\bf Global existence and uniqueness of solutions to a model of price formation}
\author{Lincoln Chayes\thanks{Department of Mathematics, University of California Los Angeles, Los Angeles, CA 90095. Email: lchayes@math.ucla.edu}, \hspace{1 pt} Mar\'ia del Mar Gonz\'alez\thanks{Universitat Polit\`ecnica de Catalunya, ETSEIB - Departament de Matem\`atica Aplicada I,  Av. Diagonal 647,
08028 Barcelona, SPAIN. Email: mar.gonzalez@upc.edu}, \\
Maria Pia Gualdani\thanks{Department of Mathematics, University of Texas at Austin, 1 University Station C1200 Austin, Texas 78712, USA. Email: gualdani@math.utexas.edu} \hspace{4 pt}and \hspace{2 pt}Inwon Kim \thanks{Department of Mathematics, University of California Los Angeles, Los Angeles, CA 90095. Email: ikim@math.ucla.edu}}
\date{}
\maketitle

\vspace{-.125 in}

\abstract{
We study a model due to J.M. Lasry and P.L. Lions
describing the evolution of a scalar {\it price}
which is realized as a free boundary in a 1D
diffusion equation with dynamically evolving, non--standard sources.  We establish global existence and uniqueness.
}

\setlength\parskip{.04 in}
\section{Introduction}

\setcounter{equation}{00}

Here we are concerned with the following PDE:
$$\left\{
\begin{array}{lll}
f_t-f_{xx} =  
\left[\delta_{p(t)+\underline{a}} - \delta_{p(t)-\underline{a}} \right]
f_x(p(t),t)&\hbox{ in }& (-1,1)\times[0,\infty);\\ \\
f_x(1,t) = f_x(-1,t)=0. && \\ \\
f(x,0)=f_I(x).&&
\end{array}\right.\leqno (\text{P})
$$
where $p(t)=\{x:f(x,t) =0\}$ presumed, for a.e.~$t$, to be a singleton, and $\underline{a}= \min\{a, |p(t)\pm 1|\}$ with $a<1$.

The model with $\underline{a} \equiv a$ was invented in \cite{Lasry-Lions} and, as explained therein
(see also \cite{price-formation},\cite{price-formationII}) is purported to describe the
dynamic evolution of a {\it price} $p(t)$ as influenced by a population of buyers and sellers.  In this initial reference, the existence of solutions was discussed, mostly in the context of a non--compact domain.

While the model on compact domains was featured (strictly speaking, the model on $\mathbb R$
does not make economic sense) there was no proviso for the circumstance $|1 \pm p(t)| < a$.
Our modification using $\underline{a}$
provides this definition and later (see Lemma~\ref{lem:stayaway}) allows us to ensure that $p(t)$ stays away from the domain boundaries at $x = \pm 1$. In terms of the model, our modification can be viewed as a ``rescue plan'' to prevent prices from severe deflation or inflation.

Recently, \cite{price-formation} the problem was solved completely for the case of symmetric initial data and in the work \cite{price-formationII}, global existence, uniqueness and stability was established for initial data sufficiently close (in a certain sense) to the piecewise linear equilibrium solution.
\hspace{-8 pt}
\footnote
{It is noted -- but not proven -- that {\it existence} might be established via a connection to a stochastic interacting particle model.  The utility of this connection is under investigation by the authors, particularly with regards the question of global stability.}
Finally, contemporaneous to the present work, a  regularized version of (P) for the non--compact case is reinvestigated in \cite{Markowich:price-formation}.  A complete derivation of uniqueness for short times -- roughly the equivalent
of our \S 1 -- is presented therein.

Notwithstanding the benign appearance of (P), the system contains intrinsic and convoluted non--linearities.  Indeed, the driving term at the sources is the gradient at the dynamically generated zero
(free boundary) of $f$ -- which in turn controls the location of the sources.
Thus, a central technical issue is to establish non--degeneracy at the free boundary and thereby
some degree of control for its motion.
E.g., in this context Hopf's lemma, while useful, is not immediately decisive without some additional regularity information at the boundary.

We consider the initial data $f_I\in C^2([-1,1])$ satisfying the following:

\begin{enumerate}
\item[(i)] $\{f_I(x)=0\} = \{p_I\},$  $f_I(x)>0$ for $x<p_I$ and $f_I(x)<0$ for $x>p_I$.
\item[(ii)] $\partial_x f_I(-1)=\partial_x f_I(+1)=0$.
\item[(iii)] Given $\lambda_I:=-\partial_x f_I(p_I)>0$, we must have
\be -\partial_x f_{I}(x) > \frac{1}{2}\lambda_{I} \quad\mbox{in }(p_I-a_0,p_I+a_0)\label{choose-a-0}\ee
for some $0<a_0 \ll a/4$.
\end{enumerate}

It is worthwhile to notice that problem (P) satisfies the following conservation identities
\begin{align}\label{mass_conserv}
\int_{-1}^{p(t)} f(x,t)\;dx = \int_{-1}^{p_I} f_I(x)\;dx=M_b>0,\quad \int_{p(t)}^{1} f(x,t)\;dx = \int_{p_I}^{1} f_I(x)\;dx= -M_p < 0.
\end{align}
Note also that the free boundary moves with a velocity given by the formula:
$$\dot{p}(t) = - \frac{f_{xx}(p(t),t)}{f_x(p(t),t)},$$
and the flux across the free boundary is given by
$$\lambda(t):=-\partial_x f(p(t),t).$$

The main result in this paper is written below:

\begin{theorem}\label{thm:main}{\rm [Global existence and uniqueness of classical solution]}
\medskip

Consider the system described on ($P$) with initial data $f_{I}$ satisfying the conditions (i)--(iii) above, and let us define  $\Omega:=-(-1,1)\times (0,\infty)$. 
Then there is a unique function $f(x,t)\in L^\infty(\Omega) $ satisfying the following:

\medskip

\noindent (A) For any $0<t\leq T$ there exists $r=r(T)$ 
with $r > 0$
such that\\
\begin{itemize}
 \item[(a)]$f$ is $C^\infty$ in $\{(x,t):|x-p(t)|\leq r\}$;\\
 \item[(b)] $\lambda(t)>r$; \\
\item[(c)] $p(t)\in (-1+r,1-r).$
\end{itemize}

\medskip

\noindent (B) $f$ solves the first two equation of $(P)$ in the classical sense (in terms of Duhamel's formula ) in $\Omega$, and $f(x,t)$ uniformly converges to $f(x,0)$ as $t\to 0$.
\end{theorem}


\section{Short times}

\setcounter{equation}{00}
The preliminary results are based on the short--time contraction principle of the following iteration:
given $f_n(x,t)$ and $p_n(t)$ such that $p_{n}(t):=\{x: f_n(\cdot,t)=0\}$ consists of a unique point for each $t>0$, consider the function $\lambda_n(t)$ defined by
$$\lambda_n(t) := - ({f_n})_x(p_n(t),t).$$
Let $f(x,t)$ solve

\begin{eqnarray}
 \pd{f}{ t} -\frac{\partial^2 f}{\partial x^2}&=&\lambda_n(t)\left[ \delta_{x=p_n(t)-\underline{a}}-\delta_{x=p_n(t)+\underline{a}}\right],\\
f_x(-1,t)&=&f_x(1,t)=0,\label{linear-problem} \\
f(x,0)&=&f_I(x).\nonumber
\end{eqnarray}
where in the above, it is tacitly assumed that for all $t$
in $(0,t_{0})$, $p_{n}(t)$ is a single point.
Then the solution to the above becomes
$f_{n+1}(x,t)$, i.e., serves to define $p_{n+1}(t)$ and $\lambda_{n+1}(t)$.

\medskip

Specifically, we define the map
\begin{eqnarray}\label{map}
\Phi :  \; L^\infty((0,t_0);X) & \to& L^\infty((0,t_0);X) , \\
 f_n &\mapsto& f_{n+1},
\nonumber
\end{eqnarray}
where $X:=  L^\infty(p_I-a_0,p_I+a_0)$ for some $a_0$ sufficiently small. 

\medskip

First let us write the solution $f(x,t)$ of \eqref{linear-problem} using the Duhamel formula:
\be\label{Duhamel}\begin{split}
f(x,t) &= \int_{-1}^{+1} \Gamma(x,x^{\prime};t)f_I(x^{\prime})\;dx^{\prime}\\& +
\int_0^t \left[\Gamma(x,(p_n(t^{\prime})-\underline{a});t-t^{\prime}) - \Gamma(x,(p_n(t^{\prime})+\underline{a});t-t^{\prime}) \right]\lambda_n(t^{\prime})
\;dt^{\prime}\\
&=:I_1+I_2,
\end{split}\ee
where $\Gamma$ denotes the fundamental solution appropriate for the domain
$(-1, +1)$ with Neumann boundary conditions:
\begin{equation}
\Gamma(x,x^{\prime}; t)  =
\sum_{k = -\infty}^{\infty}K(x - (2k + [-1]^{|k|}x^{\prime}),t)
\end{equation}
with
 $$K(x,t) = \frac{1}{\sqrt{4\pi t}}e^{- \frac{x^2}{{4 t}}}.$$

The main result of this section is stated in the following:
\begin{theorem}\label{thm-contraction}
There exists a time $t_0$  depending only on $\|f_I\|_{L^\infty(-1,1)}$ and $\lambda_I$
such that
\begin{align*}
\sup_{t\in[0,t_0]} \| f-g\|_X \le \frac{1}{2}\sup_{t\in[0,t_0]} \| f_n-g_n\|_X ,
\end{align*}
where  $X:= L^\infty(p_I-a_0,p_I+a_0)$ and the functions $f,g$ solve (\ref{linear-problem}) with right-hand side, respectively, $(f_n,p_n)$ and $(g_n,q_n)$.
\end{theorem}

\bigskip

\begin{corollary}
There exists $t_0>0$ depending only on $\|f_I\|_{L^\infty(-1,1)}$ and $\lambda_I$ such that (P) has a solution $f,p$ for all times $t\in[0,t_0]$. Moreover, $f$ is smooth accross the free boundary, and 
$$\lambda(t_0)\geq \frac{\lambda_I}{4}>0.$$
\end{corollary}

\bigskip

The proof of theorem \ref{thm-contraction} is a consequence of the following series of results:

\medskip

\begin{lemma}\label{Gamma}
For all $x$ and $y$ such that  $|x-y|<a/4$, it holds that
$$\Gamma(x,y\pm a,t) + \;|\Gamma_x(x,y\pm a,t)| + \; |\Gamma_{xx} (x,y\pm a,t)| \le G(a),$$
where the constant $G(a)$ depends only on $a$.
\end{lemma}

\begin{proof}
It follows from elementary analytical considerations. We show it for $\Gamma_x$, the rest are similar. Writing
$$\vartheta_{k}  =  \frac{x - (2k + [-1]^{|k|}x^{\prime})}{2\sqrt{t-t^{\prime}}}$$
we have
\begin{equation}
\label{fund-solut-estimate}
\Gamma_{x} (x,x';t) =  -\frac{1}{\sqrt{4\pi}}\sum_{k}\frac{1}{t-t^{\prime}}
\vartheta_{k}\text{e}^{-\vartheta_{k}^{2}}
=G_{1} \sum_{k}\frac{1}{(x - (2k + [-1]^{|k|}x^{\prime})^{2}}
\vartheta_{k}^{3}\text{e}^{-\vartheta_{k}^{2}}
\end{equation}
with $G_{1}$ a constant.
Obviously, for any $k$, we can bound $G_{1}\theta_{k}^{3}\text{e}^{-\theta_{k}^{2}}
\leq G_{2}$ for another constant $G_{2}$, whatever the value of $\vartheta_k$ might be.

\medskip

Next, since for the relevant $x^{\prime} = y \pm a$ we have that
$|x - x^{\prime}| \geq a/2$
and $|x + x^{\prime} - 2| \geq {g}_{3}a$
for a constant $g_{3}$,  we may sum the series replacing $x\pm x^{\prime}$ by the relevant worst case scenarios.
It is concluded that $|\Gamma_{x}|$ -- unintegrated in the $t$--variable --
is bounded by a finite constant (which depends on $a$, but not on $x$, $x'$ or $t$).
\end{proof}

\begin{lemma}\label{lemma_lambdak}
Consider
\begin{align*}
f_x(x,t)& = \int_{-1}^1 \Gamma_x(x,x^{\prime};t)f_I(x^{\prime})\;dx^{\prime}\\ \\
& +
\int_0^t \left[\Gamma_x(x,p_n(t^{\prime})-a;t-t^{\prime}) - \Gamma_x(x,p_n(t^{\prime})+a;t-t^{\prime}) \right]\lambda_n(t^{\prime})\;dt^{\prime}
\\ \\
&= : \frac{d I_1(x,t) } {dx}(x,t) + \frac{d I_2(x,t) } {dx}(x,t).
\end{align*}
If $\bar t$ is a time such that $|p_{I}-p_n(t)| < \bar a_0<a/4$ for all $t<\bar t$, then we can estimate
\begin{align*}
 \abs{\frac{d I_2(x,t) } {dx} } & \le G(a)\int_{0}^{t} | \lambda_n(s)|\;ds,
 \end{align*}
for all $t < \bar t$ and $x\in (p_I-\bar a_0,p_I+\bar a_0)$, where $\lambda_I>0$, $p_I$ are  as described in the hypothesis (i) -- (iii), and $G(a)$  is a constant depending only on $a$.

Moreover, for $x\in (p_I-a_0,p_I+a_0)$, the linear term can be estimated as
\begin{align*}
\frac{d I_1(x,t) } {dx} & \le -\frac{1}{2}\lambda_I e^{-t}, \quad \textrm{for all} \; t\le \sqrt{a_0}.
\end{align*}
\end{lemma}

\medskip

\begin{proof}

The second estimates is a direct consequence of
Harnack inequality. Indeed, $-\frac{dI_1}{dx}$ is a solution of an homogeneous heat equation with initial data  $-\partial_x f_I(x)$, that has a definitive sign (and lower bound) throughout $(p_{I} - a_{0}, p_{I} + a_{0})$.  Thence we can apply the classical parabolic Harnack inequality from \cite{Moser:parabolic-Harnack} to $-\frac{dI_1}{dx}$.

\medskip

The first estimate follows easily from  Lemma \ref{Gamma}, just taking into account that
$$\abs{x-p_n(t)}\leq \abs{x-p_I}+\abs{p_I-p_n(t)}\leq 2a_0 <<a.$$
\end{proof}

\bigskip

\begin{corollary}\label{cor-derivative-bound}
Suppose, in the $n^{\text{th}}$ stage of the iteration, that for all
$t < \bar t$ we have
\be\label{hypothesis2}|\lambda_{n}(t)|  \leq 2\frac{ \|f_I\|_{L^\infty(-1,1)}}{\sqrt{t}}\quad \mbox{and}\quad \abs{p_I-p_n(t)}<\bar a_0\ee
Then we can find the following bound for the next step $f=f_{n+1}$, given in \eqref{Duhamel}: there exists a time $t_0\leq \bar t$, depending on $\norm{f_{I}}_{L^{\infty}(-1,1)}$
and $\lambda_I$, such that
$$\inf_{t,x} [-f_x(x,t)] \geq \frac{\lambda_I}{4},$$
for all $t\le t_0$ and $|x-p_I|\le \min\{a_0,\bar a_0\}$.
\end{corollary}
\begin{proof}
This is an obvious consequence of the two bounds from Lemma \ref{lemma_lambdak} and the hypothesis \eqref{hypothesis2}.  Indeed we may write, uniformly for
$x \in (p_{I} - a_{0},p_{I} + a_{0})$
\begin{equation*}
(f_{n+1})_{x}(x,t)  \leq  -\frac{1}{2}\lambda_{I}\text{e}^{-t} + 4 G(a)
 \|f_I\|_{L^{\infty}(-1,1)}\sqrt{t}.
\end{equation*}
It is clear that for $t < t_{0}$ with $t_{0}$ sufficiently small, and depending only on the quotient $\norm{f_I}_{L^{\infty}(-1,1)}/\lambda_I$, the desired bound will hold.  It is noted that the time $t_0$ does not depend on the iteration coefficient $n$.
\end{proof}

In the next two lemmas we show a $L^\infty$--bound for $f_x$ and $f_{xx}$ in a neighborhood of the free boundary which is the required input for Corollary \ref{cor-derivative-bound}.
 \begin{lemma}\label{estim_fx}
 Suppose that $\bar t $ is a time such that $|p_{I}-p_n(t)| < \bar a_0$ for all $t<\bar t$ and that
\be\label{hypothesis1}\sup_x|({f_n})_x(x,t)|\le 2\frac{\|f_I\|_{L^\infty(-1,1)}}{\sqrt{t}}, \quad 0<t<\bar {t},\quad |x-p_I|< \bar a_0,\ee
 then
 $$
 \sup_{x} |(f_{n+1})_x(x,t)|\le 2\frac{\|f_I\|_{L^\infty(-1,1)}}{\sqrt{t}}, \quad 0<t<\bar {t},\quad |x-p_I|< \bar a_0 .$$
 \end{lemma}
 \begin{proof}
As before, we let
 \begin{align*}
f_x(x,t)& = \int_{-1}^1 \Gamma_x(x,x^{\prime};t)f_I(x^{\prime})\;dx^{\prime}\\& +
\int_0^t \left[\Gamma_x(x,(p_n(t^{\prime})-a);t-t^{\prime}) - \Gamma_x(x,p_n(t^{\prime})+a;t-t^{\prime}) \right]\lambda_n(t^{\prime})\;dt^{\prime}\\
&= : \frac{d I_1(x,t) } {dx} + \frac{d I_2 } {dx}(x,t).
\end{align*}
Regularity estimates for caloric functions (see \cite{DiBenedetto}, Chapter V, Theorem 8.1) imply
$$\sup_x \left|  \frac{d I_1(x,t) } {dx} \right| \le \frac{\|f_I\|_{L^\infty(-1,1)}}{\sqrt{t}}, \quad \textrm{for}\;  |x-p_I|<\bar a_0.$$
 For the nonlinear part we make use of the estimates for the Green's function given in Lemma \ref{Gamma} and the hypothesis \eqref{hypothesis1}. It holds
 \begin{align*}
\sup_x\abs{ \frac{d I_2(x,t) } {dx} } & \le G(a) \int_0^t \abs{\lambda_n(s)}\;ds \le 4 G(a)  \|f_I\|_{L^\infty(-1,1)} \sqrt{t}, \quad |x-p_I|\le \bar a_0.
\end{align*}
As soon as  $\bar t $ is small enough, say less than $\frac{1}{4G(a)}$, we get
$$\sup_x | f_x(x,t)| \le  2\frac{\|f_I\|_{L^\infty(-1,1)}}{\sqrt{t}}.$$
for all $t\leq \bar t$.
\end{proof}

The last bound we need is to show is the following local estimate for $\partial_t f$:
\begin{lemma}\label{estim_ft}
 Suppose that $\bar t $ is a time such that $|p_{I}-p_n(t)| < \bar a_0$ for all $t<\bar t$. For any $ |x-p_I|<\bar a_0$, the following holds:
$$\sup_x|\partial_t f(x,t)|\le 2\frac{\|f_I\|_{L^\infty(-1,1)}}{\sqrt{t}},\quad t<\bar t.$$
\end{lemma}

\begin{proof}
We use a very similar argument to the one of the previous lemma. Write the solution as \eqref{Duhamel}.
For the linear part we use again classical estimates
(c.f.~\cite{DiBenedetto}, chapter V, Theorem 8.1) that imply
$$|\partial_t{I_1}| \le  \frac{\|f_I\|_{L^\infty(-1,1)}}{\sqrt{t}}.$$
For the nonlinear part, consider
$$\abs{\partial_t{I_2}} =  {\partial_{xx}{I_2}}\leq G(a)\int_0^t \abs{\lambda_n(s)} \le 4 G(a){\|f_I\|_{L^\infty(-1,1)}}{\sqrt{t}},$$
where we have used Lemma \ref{Gamma} for the estimate of $\Gamma_{xx}$ and \eqref{hypothesis1}.\\
Choosing $\bar t$ as in the previous lemma, it follows that
$$
\sup_x|\partial_t f(x,t)|\le 2\frac{\|f_I\|_{L^\infty(-1,1)}}{\sqrt{t}} \quad \mbox{for all }t\leq \bar t.
$$
\end{proof}

\noindent\emph{Remark: }
Note that the time $\bar t$ depends only on some universal constants, and the time the free boundary uses to exit an interval. Unfortunately in the proof of Theorem \ref{thm-contraction} we will need to choose a smaller time step $t_0< \bar t$, which depends on the initial data.  

\medskip

The above lemma gives us an estimate in terms of the $L^\infty$--norm of the initial data. \\
Is is also possible to estimate the $f_x$ in terms of the $L^1$--norm (which is in our case constant in time), at the price of a worse denominator. We will not make use of this lemma in this section, but later in the proof of global existence:

\begin{lemma}\label{estim_fx_L1}
 Suppose that $\bar t $ is a time such that $|p_{I}-p_n(t)| < a/4$ for all $t<\bar t$ and that
\be\label{hypothesis1}\sup_x|({f_n})_x(x,t)|\le 2\frac{\|f_I\|_{L^1(-1,1)}}{(\sqrt{t})^3}, \quad 0<t<\bar {t},\quad |x-p_I|< a/4,\ee
 then
 \begin{align*}
 \sup_{x} |(f_{n+1})_x(x,t)|\le 2\frac{\|f_I\|_{L^1(-1,1)}}{(\sqrt{t})^4}, \quad 0<t<\bar {t},\quad |x-p_I|< a/4 ,\\
  \sup_{x} |(f_{n+1})_{xx}(x,t)|\le 2\frac{\|f_I\|_{L^1(-1,1)}}{(\sqrt{t})^5}, \quad 0<t<\bar {t},\quad |x-p_I|< a/4.
 \end{align*}
 \end{lemma}
 \begin{proof}
 The above estimates can be proven following the same steps as in Lemma \ref{estim_fx} and \ref{estim_ft}. 
Instead of using the regularity estimates for the caloric function with the $L^\infty$--norm of the function, we make use now of the $L^1$--norm. From  \cite{Evans} (Theorem 9, section 2.3.c.) we know that
\begin{align*}
\sup_x \left| \frac{dI_1}{dx} \right|  \le \frac { \|f\|_{L^1}}{t^2},\\
\sup_x \left| \frac{d^2I_1}{dx^2} \right|  \le \frac { \|f\|_{L^1}}{t^{5/2}},
\end{align*}
where $I_1$ is defined as in Lemma  \ref{estim_fx}. The rest of the proof follows similarly. \\

 \end{proof}

In the next corollary we show that given $p_n(t)\in (p_I-a_0/2, p_I+a_0/2)$, also the free boundary $p(t)$ of the next iteration step is well defined and moreover it is $C^{1/2}$--H\"{o}lder continuous. In particular, also $p(t)\in  (p_I-a_0/2, p_I+a_0/2)$ for $t$ small enough.

\begin{corollary}\label{corollary-velocity}
There exists a time $t_0>0$ such that the solution $f(\cdot,t)$ constructed in \eqref{Duhamel} has a unique zero in $(p_I-a_0,p_I+a_0)$, for each fixed time $t<t_0$, which we denote by $p(t)$.  Moreover, $p(t)$ satisfies
$$|p(t)-p_I | \le a_1\sqrt{t},$$
where $a_1$ only depends on the initial data as shown in the proof.
\end{corollary}

\begin{proof}
First, note that Corollary \ref{cor-derivative-bound}
implies monotonicity of $f(\cdot, t)$ in $(p_I-a_0,p_I+a_0)$ which
assures the existence of at most one zero, $p(t)$. Next,  since $f(p(t),t)=0$ for all $t$, then
\begin{equation}
\label{ODE}
\dot{p}(t) =-\frac{f_t(p(t),t)}{f_x(p(t),t)}.
\end{equation}
On the other hand, Corollary \ref{corollary-velocity} gives a bound for the velocity $|f_t(x,t)|\le2 \frac{\|f_I\|_{L^\infty(-1,1)}}{\sqrt{t}}$ in the time interval $(0,\bar t)$, and Corollary \ref{cor-derivative-bound} bounds the slope $\abs{f_x}$ from below in an interval $\abs{x-p_I}<a_0$ in $(0,t_0)$. Existence results for the ODE in Eq. (\ref{ODE})  implies the existence of this $p(t)$.

As a consequence, the interval $|p(t)-p_I |$ is bounded by
$$|p(t)-p_I | \le  8\frac{\|f_I\|_{L^\infty(-1,1)}}{\lambda_I}{\sqrt{t}},\quad \textrm{for} \; 0<t\le t_0.$$
The thesis follows by choosing the constant $a_1=8\frac{\|f_I\|_{L^\infty(-1,1)}}{\lambda_I}$, and the time $t_0$ as
\be\label{choose-t-0} 8\frac{\|f_I\|_{L^\infty(-1,1)}}{\lambda_I}\sqrt{t_0} =: \frac{a_0}{2}.\ee
\end{proof}

\medskip

 Now we are ready to show that the map $\Phi$ defined in Eq.\eqref{map} is a contraction in the space $L^\infty\lp(0,t_0);X\rp$, and consequently to define $t_0$. At this point we already have a specific condition on the time $t_0$: it must satisfy all the conditions from the previous results, of which the most restrictive is precisely Eq.\eqref{choose-t-0}.\\

\bigskip

\noindent\textbf{Proof of Theorem \ref{thm-contraction}: }

  First note that if $f_n\in L^\infty([0,t_0], X)$ with $f_n(p_n)=0$  satisfies
  \be\label{conditions-iteration}p_n(t)\in(p_I-a_0/2,p_I+a_0/2), \quad \sup_{x} | {f_n}_x(x,t)|\le 2 \frac{\|f_I\|_{L^\infty(-1,1)}}{\sqrt{t}}, \quad {f_n}_x(p_n,t) \le -\frac{\lambda_I}{4},\ee
then the previous results assure that the next step in the iteration $f(x,t)$ has a well defined zero $p(t)$ and satisfies also the same estimates.

\medskip

Let $f_n, \; g_n \in  L^\infty([0,t_0], X)$ be such that $f_n(p_n)=0$, $g_n(q_n)=0$, $\lambda_n=-f_n(p_n)$, $\xi_n=-g_n(q_n)$, and that both satisfy
Eq.\eqref{conditions-iteration}. Then we estimate the difference of the images by our mapping $\Phi$: for any $x\in(p_I-a_0,p_I+a_0)$,
\begin{align*}
 f(x)-g(x) = & \int_0^t (\Gamma(x,p_n-a;t-s)-\Gamma(x,p_n+a; t-s))\lambda_n(s)\;ds \\
 &-  \int_0^t (\Gamma(x,q_n-a;t-s)-\Gamma(x,q_n+a; t-s))\xi_n(s)\;ds \\
 = &  \int_0^t \left[\Gamma(x,p_n-a;t-s)-\Gamma(x,p_n+a;t-s)\right](\lambda_n(s) - \xi_n(s)) \;ds \\
 & +   \int_0^t \xi_n(s)[\Gamma(x,p_n-a;t-s)-\Gamma(x,q_n-a; t-s) \\
 &\quad\quad\quad\quad+\Gamma(x,p_n+a;t-s)-\Gamma(x,q_n+a; t-s)]  \;ds\\
 =&: A_1+A_2.
\end{align*}
It holds
\begin{align*}
\abs{A_1}\le& \sup(t) | \lambda_n(t)-\xi_n(t)| \int_0^t | \Gamma(x,p_n-a;t-s)-\Gamma(x,p_n+a; t-s) |\;ds.
\end{align*}
We now need to estimate the difference $ | \lambda_n(t)-\xi_n(t)|$ in terms of the quantity $\| f_n-g_n \|_{L^\infty(p_I-a_0, p_I+a_0)}$. Therefore
\begin{align*}
|\lambda_n(t) - \xi_n(t)| \le &\abs{{f_n}_x(p_n(t),t) - {g_n}_x(\xi_n(t),t)} \\ \\
\le &|  {f_n}_x(p_n(t),t) - {g_n}_x(p_n(t),t) | +  |  {g_n}_x(p_n(t),t) - {g_n}_x(\xi_n(t),t) | \\ \\
\le & \| f_n-g_n\|_{C^1 (p_I-a_0,p_I+a_0)} + \|{g_n}_{xx}\|_{ L^\infty (p_I-a_0,p_I+a_0)}|p_n(t) - q_n(t)|.
\end{align*}

\medskip

Next we estimate the term $ \|{g_n}_{xx}\|_{ L^\infty (p_I-a_0,p_I+a_0)}$ with Lemma \ref{estim_ft}.
and the $C^2$--norm of $f_n-g_n$ with classical estimates for the caloric function. Recall that the estimates hold in a neighborhood of the free boundary far away from the source and sink.  Therefore
\begin{align*}
|\lambda_n(t) - \xi_n(t)| \le &  \frac{\| f_n-g_n\|_{L^\infty (p_I-a_0,p_I+a_0)} }{\sqrt{t}}+  2 \frac{\|f_I\|_{L^\infty(-1,1)}}{\sqrt{t}} |p_n(t) - q_n(t)|,
\end{align*}
for each fixed time $t$.
On the other hand, by an elementary geometrical argument, it holds
\begin{align*}
|p_n(t) - q_n(t)| \le \sup_x | f_n(x,t) - g_n(x,t) | \cdot \frac{4}{\lambda_I},
\end{align*}
since the slope is bounded by $\frac{\lambda_I }{4} \le | {f_n}_x(x,t)| \le 2 \frac{\|f_I\|_{L^\infty(-1,1)}}{\sqrt{t}}$.\\
Hence
\be\label{formula10}\abs{\lambda_n(t) - \xi_n(t)}\leq  \left[ 1 +\frac{8}{\lambda_I} \|f_I\|_{L^\infty(-1,1)} \right]\frac{
 \| f_n-g_n\|_{L^\infty (p_I-a_0,p_I+a_0)}}{\sqrt{t}}. \ee
Consequently, the term $A_1$ can be estimated using \eqref{formula10} and Lemma \ref{Gamma} as
\begin{align*}
 \abs{A_1} \le 2G(a) \left[ 1 +\frac{8}{\lambda_I} \|f_I\|_{L^\infty(-1,1)} \right]{
 \| f_n-g_n\|_{L^\infty (p_I-a_0,p_I+a_0)}}{\sqrt{t}}.
\end{align*}

 Similar arguments can be applied in the estimates for $A_2$: the hypothesis \eqref{conditions-iteration} gives that
 \begin{align*}
\abs{A_2} &\le  2 \frac{\|f_I\|_{L^\infty(-1,1)}}{\sqrt{t}} \int_0^t |\Gamma(x,p_n-a;t-s)-\Gamma(x,q_n-a; t-s))|\;ds \\
 &+  2 \frac{\|f_I\|_{L^\infty(-1,1)}}{\sqrt{t}} \int_0^t |\Gamma(x,p_n+a;t-s)-\Gamma(x,q_n+a; t-s))|\;ds
 \end{align*}
 Next, note that $(p_n-a, q_n-a) \cap  (p_I-a_0,p_I+a_0) = \varnothing$. Here we are slightly abusing the notation by assuming that $p_n<q_n$. It holds, for each $t<t_0$, that
\begin{align*}
\int_0^t &|\Gamma(x,p_n-a;t-s)-\Gamma(x,q_n-a; t-s))|\;ds \\
&\le |p_n - q_n| \sup_{ y\in (p_n-a, q_n-a)} \int_0^t |\Gamma_y(x,y; t-s)|ds\\
&\le G(a) \frac{4}{\lambda_I} \| f_n-g_n\|_{L^\infty (p_I-a_0,p_I+a_0)}t,
\end{align*}

 In conclusion we have
 $$ |f(x)-g(x)| \le c \sqrt{t} \| f_n-g_n\|_{L^\infty(p_I-a_0,p_I+a_0)},$$
where the constant $c$ is given by
$$c:= G(a) \left[ 1 +\frac{8}{\lambda_I} \|f_I\|_{L^\infty(-1,1)} \right].$$
\noindent Choosing $t_0$ small enough such that $$G(a) \left[ 1 +\frac{8}{\lambda_I} \|f_I\|_{L^\infty(-1,1)} \right]\sqrt{t_0}<1,
$$
the proposition is proved.
\qed\\

\noindent \emph{Remark: }  We remark
that the proof in \cite{price-formationII} for short time existence -- almost exactly the approach of
\cite{Markowich:price-formation} -- made use of semi--group methods.  The derivation here, while admittedly less sophisticated, is certainly more robust and is deliberately tailored to the upcoming developments.  In particular, we
now have tangible criteria under which short time existence is purported to break down.  Indeed, a parallel derivation for the regularized problems allows, at least for the compact case, an immediate proof of global existence:  With the asset of regularization, flux bounds easily imply various derivative bounds which are the central objective of the next sections.


\section{An $L^1$ bound on the flux (and the rescue plan at the Neumann boundary)}

In this section we show that, as long as the solution of (P) exists, the flux of the solution at the zero set stays bounded. This  result will be then used in the next section to provide further estimates on the derivatives of the solutions.  As a corollary, we will also show that the zero set of the solution cannot approach  the Neumann boundaries too closely.

To facilitate matters, we shall, in essence decouple the positive and negative pieces of $f$ and, in addition, describe problems of this sort on a larger space which restores much of the linearity usually associated with diffusion problems.  

Thus, first, we shall define (P$^{\prime}$) to be a one--sided version of the system (P) that is to say (i)  The positive part of $f$ is set, identically, to zero (and no source).  (ii)  The zero $p(t)$ is predetermined.  Thus, in principle, in order to recover the system (P), two such (P$^{\prime}$) models can be glued together subject to ``additional constraints'' on their mutual $p(t)$.  For the single sided problem, usually to the right of $p(t)$,
we will use
the notation 
$\rho_{p}(x,t)$ for the density, $M_{p}$ for the total mass, etc.   
When two such models are to be used in tandem, the one on the left (now representing the {\it positive} part of $f$) will be denoted by $\rho_{b}$.

Secondly, we may define these sorts of systems -- (P$^{\prime}$) will be sufficient -- on a 
{\it foliated} space.  Let $\mathbb N$ denote the natural numbers, including zero, and consider 
$[-1,+1]^{\mathbb N}$.  We shall refer to the individual 
elements as {\it levels} and, denote these, along with various associated quantities with a superscript:
$[-1,+1]^{(0)},  [-1,+1]^{(1)}, \dots $  
Let $p(t): [0,T]  \to (-1,+1)$ denote a continuous function.
On each level, we have a copy of $p(t)$ (always ``located'' in the
corresponding position) and we consider a sequence of densities
$\rho_{p}^{(0)}(x, t), \ \rho_{p}^{(1)}(x,t), \dots, $ 
with the $n^{\text{th}}$ density supported in
$[p(t), 1]^{(n)} \subset [-1,1]^{(n)}$.  Initially,
$\rho_{p}^{(n)}(x, 0) \equiv 0$ for $n > 0$ while
$\rho_{p}^{(0)}(x, 0)  =  \rho_{p}(x, 0)$.  Each of the
$\rho_{p}^{(n)}(x, t)$'s obey the diffusion equation with a source to be described below,
Dirichlet boundary conditions at their respective $p(t)$ and Neumann condition at the corresponding $x = 1$.  Finally, each $\rho_{p}^{(n)}(x, t)$, $n > 1$ has a source which is located at its respective
$p(t) + \underline{a}$ (where,
we remind the reader,
$\underline{a} = \min\{a,\frac{1}{2}[1 - p(t)] \}$)
and has strength
provided by $\lambda^{(n)}(t)  =  \nabla\rho_{p}^{(n-1)}(p(t), t)$.


In the context of the foliated model, we may calculate various quantities for the original model.
Of particular relevance, it is seen that:
\begin{equation}
\label{BSE}
\sum_{n = 0}^{\infty}\rho_{p}^{(n)}(x, t)  =  \rho_{p}(x, t)
\end{equation}
and
\begin{equation}
\label{ASE}
\int_{0}^{t}\lambda(t^{\prime})dt^{\prime}  =
\sum_{n = 0}^{\infty}n
\int_{p(t)}^{1}
\rho_{p}^{(n)}(x, t)dx.
\end{equation}

To vindicate the above claims, it is first necessary to demonstrate that for $t \in [0,T^{\prime}]$ with
$T^{\prime} < T$, the tail of 
$[-1,+1]^{\mathbb N}$ is sparsely populated.  This is elementary and follows {\it mutatis mutandis} from arguments in the proof of Proposition \ref{VV} 
where detailed estimates along these lines are provided.
Then, defining
\begin{equation}
\tilde \rho_{p}^{[N]}  =  \sum_{n = 0}^{N}\rho_{p}^{(n)}
\end{equation}
we have, by monotonicity,  
$$
\tilde \rho_{p} \lim_{N\to \infty} \tilde \rho_{p}^{[N]}
$$
satisfying the system 
(P$^{\prime}$) which demonstrates Eq.(\ref{BSE}).


As for Eq.(\ref{ASE}), we define
\begin{equation}
M_{p}^{(n)}(t)  = \int_{p(t)}^{1}\rho_{p}^{(n)}(x,t)dx
\end{equation}
it is seen that $\sum_{n}M_{p}^{(n)}(t)$ is conserved -- and hence identically equal to $M_{p}$.  
Moreover, 
$M_{p} - M_{p}^{(0)}(t)$ is given by
\begin{equation}
M^{(0)}_{p}(0) - M^{(0)}_{p}(t)  =  
\sum_{n \geq 1}M^{(n)}_{p}(t)
=  \int_{0}^{t}\lambda^{(1)}(t^{\prime})dt^{\prime}
\end{equation}
i.e., the mass lost on the first level is exactly that which fluxed up to the higher levels.  Similarly we have for every level, the identity
\begin{equation}
\label{STK}
\sum_{n \geq \ell}M_{p}^{(\ell)}  =  
\int_{0}^{t}\lambda^{(\ell)}(t^{\prime})dt^{\prime}
\end{equation}
and summing both sides of Eq.(\ref{STK}), we obtain 
Eq.(\ref{ASE})


\begin {proposition}
\label{VV}
  For any $T < \infty$, as long as the solution of $(P)$ exists in $[0,T)$, then
\label{firstbound}
$$
\int_{0}^{T}
\hspace{-4 pt}
\lambda(t)dt < \infty.
$$
\end {proposition}
\begin{proof}
As we shall see, the consequences of a finite time divergence violate sensible notions of the slow scale for the diffusive transport of substantial material over large distances.  
The pertinent observation is that for $\delta t$ small, there is essentially no diffusion over any appreciable distance in the allotted time and in its absence, the ``essential'' supports of the $p$-- and $b$--densities become so widely separated that diffusion cannot account for
complimentary transports.

Suppose then that the above display does {\it not} hold.  Then for any
$\delta t$ and $K$ ($\delta t \ll1$ and $K \gg 1$ to be specified when necessary) there is a
$\mathfrak t_{1}$ with $\mathfrak t_{1} + \delta t < T$ such that
$$
\int_{\mathfrak t_{1}}^{\mathfrak t_{1} + \delta t}
\hspace{-8 pt}
\lambda(t)dt > K.
$$
For simplicity let us reset $\mathfrak t _{1}$ to zero and work with $t \in [0, \delta t]$.
The midrange objective is to show that under the stipulation of a large flux in a short time, $p(t)$ must head towards the boundaries.
We consider, for the time being, the one--sided perspective.
%



In this context it is noted that the foliated problem in the absence of diffusion is particularly simple: mass on each level simply gets displaced,
as $p(t)$ sweeps through it, to the next level at the position corresponding to
$p(t) + \underline{a}$.  
with the density, on each level supported to the right of the furthest excursion of $p$.
Thus, in particular if 
$a \ll 1$, and $p(t)$ does {\it not} enter the right rescue zone 
(namely $x > 1-2a$) there are at most the order of $1/a$ levels that get occupied in the sweep.  To treat the general case, let us consider for $\Delta \in (0,a)$ small but of order unity 
the system $(\text{P}_{\Delta}^{\prime})$ which is
defined exactly as the system 
$(\text{P}^{\prime})$ but with $a$ replaced by $\Delta$.  
(Also, with corresponding rescue plan, but this shall not enter into our considerations.)  Let us consider, for identical initial conditions and identical $p(t)$ which remains {\it outside the} $(\text{P}_{\Delta}^{\prime})$ {\it rescue zone}, a comparison between the 
$(\text{P}_{\Delta}^{\prime})$ and the $(\text{P}^{\prime})$ systems.  We make two claims both of which are straightforward to verify on the basis of an underlying particle model based on non--interacting stochastic elements.  The first claim is that the 
 is that the total flux in the $\Delta$--system is not smaller than that of the usual system.  The second claim is that in two such 
 $(\text{P}_{\Delta}^{\prime})$ systems with differing zeros:
 $p(t)  \leq  q(t)$ (with both staying out of the rescue zone)
 that the flux in the $p$--version is not larger than the flux in the 
 $q$--version.   
The proof of these claims follow from an obvious variant of Lemma 3.2 in \cite{Chayes-Kim:two-sided}.  

Now let us suppose that $p(t)$ -- in a (P) or ($\text{P}^{\prime}$) model -- is such that on $[0,\delta t]$, (with $\delta t$ small to be specified later)
\begin{equation}
\max_{t}p(t) < 1-4\Delta.
\end{equation}
Then we claim that the flux is limited (e.g., by the order of $\Delta^{-1})$.
We use a comparison a $\text{P}_{\Delta}^{\prime}$
system with its own $q(t)$.  Here 
$q(0)=p(0)$
and then $q$ (almost) immediately 
jets out to $x = 1 - 2\Delta$, then backs to $1 - 4\Delta$ and stays there till
time $\delta t$.
Let $N_{\Delta} \leq 2/\Delta$ denote the number of boxes in the initial surge which get occupied.  Then, by Eq.(\ref{ASE}), the initial portion of the flux is not more than $N_{\Delta}^{2}M_{p}$

Notice that (still in the $q$--system with parameter $\Delta$)
all the mass on various levels is trapped in the respective regions $\{x \geq 1 - 2\Delta\}$ with the zero a distance $2\Delta$ away.
Thus, in the remaining time, the remaining flux is determined by how much can flux across various neighborhoods of size $2\Delta$ in time $\delta t$.
To this end let us now fix $\delta t$ ``small enough''.  To be specific,
$(\delta t)^{1/4} = \gamma\Delta$ for $\gamma$ sufficiently small ensures that
the amplitude for diffusion across $2\Delta$ -- a unit mass at $1 - 2\Delta$ and a
Dirichlet zero at $1 - 4\Delta$ is less than 
$$
\varepsilon:=\varepsilon(\delta t)=\text{e}^{-1/[\delta t]^{1/2}}.
$$
We may also interpret $\varepsilon$ probabilistically:  Starting at the $k^{\text{th}}$ level, a fraction (less than) $\varepsilon$ of the initial mass makes it to the
$k + 1^{\text{st}}$ level in time $\delta t$ and a fraction (more than) $1 - \varepsilon$ stays behind.
Of the former, a fraction (less than) $\varepsilon$ gets promoted to the
$k + 2^{\text{nd}}$ level etc.
In short, the upward distribution of diffused mass after time $\delta t$ is dominated by the a geometric random variable, 
X with parameter $\varepsilon$ and so the total flux
coming from a unit mass which at $t = 0^{+}$ was on the $k^{\text{th}}$ level in $[1-2\Delta, 1]$ is bounded by
\begin{equation}
\mathbb E(\text{X})  =  \varepsilon(1 - \varepsilon) +
2\varepsilon^{2}(1 - \varepsilon) + \dots = \frac{\varepsilon}{1-\varepsilon}.
\end{equation}
This latter quantity -- independent of $k$ -- must be multiplied by the mass that was on the $k^{\text{th}}$ level at $t = 0^{+}$ and summed. 
Since the total mass on all levels adds up to $M_{p}$ this provides
a bound on the diffusive contribution
$\int \hspace {-2 pt}\lambda dt$
that is given by
\begin{equation}
\label{CEC}
\frac{\varepsilon M_{p}}{1-\varepsilon}.
\end{equation}
Thus, back in the real problem for fixed $\Delta > 0$ sufficiently small (i.e., compared with the minimum of $a$ and $(1-a)$ as will be clear below) we may choose $\delta t$ small enough and take the content of
Eq.\eqref{CEC} and add this to our $M_{p}N_{\Delta}^{2}$, and then stipulate the sum of these two to be less than one half of the total $K$.  Under these circumstances, there must be a first time $[\delta t]_{p} < \delta t$ such that
$p([\delta t]_{p}) > 1 - 2\Delta$ and moreover, at this time $[\delta t]_{p}$, the total flux is less than $\frac{1}{2}K$.

With the above in hand, we incorporate the constraints of the {\it full} problem -- namely a $\rho_{p}$ and a $\rho_{b}$ both in play with fluxes at $p(t)$ that are supposed to match. 
The first implication is that there is another such time
$[\delta t]_{b}$ at which all of $\rho_{b}$ is in $[-1, -1 + 2\Delta]$ and moreover, even at the later of these two times, there is more than half the fluxing left to be done.
Let $I_{b}$ denote the set
\begin{equation}
I_{b}  =  \{t\in[[\delta t] _{b}, \delta t] \mid p(t) < -1 + 4\Delta \}
\end{equation}
and similarly for $I_{p}$.  
(Here, finally, we choose $\Delta$ small enough so that the intervals of size $4\Delta$ about $x = \pm 1$ do not intersect)
It is clear, from the $b$--perspective that
$\int_{I_{b}^{c}}\lambda dt$ is small.  Specifically, an estimate of the form in Eq.(\ref{CEC}) applies.  Indeed, at the earlier time all of $\rho_b$ was confined to the region $x<-1+2\Delta$ and has less time, for $t\in I_b^c$, than $\delta t$ to achieve diffusion across the gap of length scale $2\Delta$.
Similarly, from the $p$--perspective the flux during $I_{p}^{c}$ is small.
Since
$I_{b} \subset I_{p}^{c}$ or vice versa,
it is seen that there is no time after
$\max\{ [\delta t]_{p}, [\delta t]_{b} \}$
in which the requisite remaining flux can be achieved.
\end{proof}
As an immediate consequence, we may conclude that -- at least in any finite time interval -- $p(t)$ stays away from a neighborhood of the Neumann boundary.

\begin{corollary}\label{lem:stayaway}
Let $T < \infty$ denote any time up to which the solution to (P) exists. Then there is an $\e>0$ (depending on $T$) such that
$$
p(t) \in [-1+\e, 1-\e] \quad\hbox{ for } 0\leq t\leq T.
$$

\end{corollary}
\begin{proof}
Suppose not. We may then assume that $p(t)$ approaches $-1$ as $t\to T$. In particular, there exists a sequence $t_n\to T$ where
\begin{itemize}
\item[(i)] $p (t_n) = \min_{[t_0,t_n+\delta_{n}]} p(t) $ for some $\delta_{n} > 0$ and  \\
\item[(ii)] $\e_n=p(t_n)+1$ satisfies $\e_{n+1} \leq \frac{1}{4}\e_n$.\\
\end{itemize}

\noindent This means that the new excursion of $p(t)$ went past the previous excursion of the source.

\medskip

Let $g_{n-1}:=\int_{t_{n-1}}^{t_n} \lambda_1(t) dt$, where $\lambda_1(t)$ is the flux from the first term in the foliation model introduced above, with initial data $f(x,t_{n-1})$.
In the context of the {\it frozen} model (ignoring diffusion effects), all the mass between $p(t_n)$ and $p(t_{n-1})$ is fluxed through the zero set at least once. On the other hand the removed mass via the flux at the zero set is then deposited onto the left side by the source term. Therefore we have
$$
g_n = \int_{[p(t_n),p(t_{n-1})]} \rho_b(x,t_{n-1}) dx \geq g_{n-1}
$$
which would be an obvious violation of Lemma~\ref{firstbound}.

\medskip

We shall show, more or less, that the diffusion cannot alter this situation in a significant way. Indeed we claim that
either $g_n \geq g_{n-1}$ or $g_n+g_{n-1} \geq c_0$, where $c_0$ is independent of $n$. Thus, at least one case occurs for infinitely many $n$, which again violates Lemma~\ref{firstbound}.

\medskip

Let us consider the situation at time $t_{n-1}$: a fraction of the mass $\rho_b$ has been fluxed and now resides in the interval $I=[-1+\frac{1}{4}\e_{n-1}, -1+\e_{n-1}]$. Let $\tau$ be the time between
$t_{n-1}$ and $t_n$ where $p(\tau) + 1 =\frac{1}{4}\e_{n-1}$. From now on, for simplicity, we denote $\e:=\e_{n-1}$.

\medskip

Note that
$$
\begin{array}{lll}
g_n \geq \int_{t_{n-1}}^{\tau} \lambda(t) dt&=&\int_{[-1+\e/4,-1+\e]} f(x,t_n) dx +\int_{t_{n-1}}^{\tau} -f_x(-1+\e/4,t) dt \\ \\
&\geq&  g_{n-1} + \int_{[-1+\e/4,-1+\e]}f(x,t_{n-2}) dx+  \int_{t_{n-2}}^{\tau} -f_x(-1+\e/4), t) dt \\ \\
&=& g_{n-1} +  \int_{[-1+\e/4,-1+\e]}f(x,t_{n-2}) dx  + A-B,
\end{array}
$$
where $A$ is the flux of $f$ diffused into $I$ from the region $-1\leq x\leq -1+\e/4$, and $B$ is the flux which diffused out from  $I$ into $x<-1+\e/4$.

 \bigskip

{\it Case 1.} First observe that if $ t_n - t_{n-2} \leq C_1\e^2$ then there is not much diffusion between $[-1+\e/4, -1+\e]$ over the time interval $[t_{n-2}, t_n]$ such that
$$
B \leq \int_{[-1+\e/4,-1+\e]} f(x,t_{n-2}) dx \hbox{ and }g_n \geq g_{n-1}.
$$

\bigskip

{\it Case 2.} Next suppose that $t_n-t_{n-2} \geq C_1\e^2$. Our goal in this case is to show that either $A>B$ or $g_n + g_{n-1} \geq c_0$.

\medskip

Let $a_n:=\int_{-1+\e/2}^{1} \rho_b(x,t_{n-2}) dx$ . Observe that, by conservation of mass,  $M_b-a_n = \int_{-1}^{-1+\e/2} \rho_b(x,t_{n-2}) dx$.  It is then not hard to see that
$$
A-B \geq  \int_{t_n}^{\tau} -h_x(-1+\e/4,t) dx
$$

where $h(x,t_n) =a_n\delta_{-1+\e/2} + (M_b-a_n) \delta_{-1}$ and $h$ solves heat equation in $[-1,1]\times [t_{n-2},\tau]$ with Neumann boundary conditions at $|x|=1$.

\medskip

The effects from $x=1$ are not significant and may be neglected for the time being. The condition at $x=-1$ is, of course, important: the equivalent situation (via reflection) is solving heat equation in $\mathbb R\times [t_{n-2},\tau]$ with initial data
$$
\tilde{h}(x,t_n)=a_n\delta_{-1+\e/2} + a_n \delta_{-1-\e/2} +2(M_b-a_n) \delta_{-1}.
$$
If we assume that $a_n\leq M_b/2$, then the flux from the sources at $x=-1$ and $x=-1+\e/2$ regulate each other, yielding
$$
h_x(-1+\e/4,t)<0\hbox{ for } t>t_n.
$$
 Therefore $A>B$.
On the other hand if $a_n\geq M_b/2$, then the flux to the $x>-1+\e/2$ is already significant so that we have $g_n+g_{n-1}$ bigger than $c_0$.  The desired result is established.

\end{proof}

\bigskip

\section{Global-time existence of the solution }

\setcounter{equation}{00}

In the previous section we showed that the problem has unique solution for a small time interval $t_0$. In this section we will show that we can iterate this process to produce the unique solution of our problem
for global times.

\medskip

Let us restart the process at the time $t_0$. This will lead to the existence (and uniqueness) of a solution in the time interval $(t_0, t_1)$; we shall continue to iterate the process as long as we can.

\medskip

The length of the time interval $t_{n+1}-t_n$ in which the contraction can take place is a function that depends on
$$t_{n+1}-t_n=C \frac{\norm{f(\cdot,t_n)}_{L^\infty(-1,1)}}{\lambda(t_n)}.$$

\medskip

Suppose there exists a {\it blow-up time} $0<t^\star<\infty$: it means that we can no longer find a small time interval during which Theorem  \ref{thm-contraction} holds.
This happens exactly when, for any sequence $t_n\to t^\star$, one of the following holds:
\begin{itemize}
\item[(i)] $\limsup_{n\to\infty}\norm{f(\cdot,t_n)}_{L^{\infty}(-1,1)}=\infty$ \quad or\\
\item[(ii)]  $\liminf_{n\to\infty}\lambda(t_n)\to 0$,\quad  or\\
\item[(iii)]  $\limsup_{n\to\infty} |f_{xx}|(p(t_n),t_n)\to\infty$.
\end{itemize}

\medskip

At the time $t^\star$ two possible configurations may happen: the limit of the free boundary $p(t)$ is an unique point as $t\to t^\star$, or a puddle of zero forms at $t=t^\star$. In both cases we will show that the arguments lead to a contradiction, showing the non-existence of such a $t^\star$.  In the non-puddle case (section 3.1) a contradiction will be yielded by proving that all relevant norms given in (i)-(iii) are bounded up to $t=t^\star$.  In the puddle case (section 3.2) the strategy is a bit different. We first show that (essentially) all derivatives of $f$ are bounded up to $t=t^\star$: we then will show that such regularity result is too strong to hold at a blow-up time, therefore concluding that the blow-up time $t=t^\star$ does not exist.

\medskip

We start showing what happens when the free boundary has an unique limit.


\subsection {Non-puddle case: $\lim_{t\to t^\star} p(t) = p(t^\star)$}
In this case we can show that, due to classical estimates on caloric functions, none of the three breaking factors listed above takes place, therefore yielding a contradiction.

\medskip

Let $t^\star$ the breaking time. When the limit is unique, 
the source $p(t)-\underline{a}$ and the sink $p(t)+\underline{a}$ also have an unique limit as $t\to t^\star$.
Moreover, the locations of the source and sink at the time $t=t^\star$ are $\underline{a}$--away from the point $x=p(t^\star)$.

\medskip

This allows to find a spacial neighborhood of $p(t^\star)$ and a time interval (before $t^\star$) in which the source and sink never entered in that neighborhood during that time interval. 
 
 \medskip
 
 Let 
 $$
 \delta:= \frac{\underline{a}(p(t^{\star}))}{8}>0
 $$ and define $I_\delta$ as a small neighborhood with radius $\delta$ of the point $p(t^\star)$, i.e., $I_\delta :=( p(t^\star)-\delta, p(t^\star)+\delta)$. Then there exists $0<t_{\delta}<t^\star$ such that 
 $I_\delta\times[t_\delta,t^\star]$ is far from any sink or source.

\medskip

In terms of estimates, the fact that sink and source are far away from $I_\delta$ implies that we can use classical estimates for the homogeneous heat equation,
 as in the proofs Lemma \ref{estim_ft} and \ref{estim_fx}.

\medskip

Let us start the problem (P) with initial datum $f_I := f(x,t_\delta)$. Following the same contraction argument for the local existence, we can show that there exists an unique solution $f(x,t)$ in the time interval $(t_\delta, t_\delta+\Delta t_1)$. Iterating the process another time, we can extend the solution to a bigger time interval  $(t^\star-\delta, t^\star-\delta+\Delta t_1+\Delta t_2)$, with $\Delta t_2 > 0$. The extension process could continue by iteration. The following is due to local estimates for solutions of the heat equation:

\begin{lemma}
Let $f(x,t)$ solve (P) in the time interval $(t_\delta, t^\star-\delta+\Delta )$. Then
\begin{align*}
|f^{(k)}(x,t^\star-\delta+\Delta)| \le C\frac {\|f(x,\cdot)\|_{L^1}}{\Delta}^{k+1/2}\hbox{ in } I_\delta.\\
\end{align*}
where $f^{(k)}$ denotes the $k$-th spatial derivative of $f$.
\end{lemma}

Due to above lemma, the mass conservation, and the iteration process, we obtain 
$$
|f^{(k)}(x,t^\star-\delta+\tau)| \le (M_p+M_b) (\tau)^{-(k+1/2)}\quad \hbox{ for } 0<\tau<t^\star - t_\delta.
$$

 In particular $\lambda(t)$ is bounded, and therefore the $f$ is uniformly bounded up to $t=t^\star$.
Therefore case  (i) and (iii) are eliminated.\\

Moreover $f_t$ is also uniformly bounded near $t=t^\star$ in $I_\delta$, and thus $f(x,t)$ uniformly converges to a $C^\infty$ function $F(x)$ as $t\to t^*$. Let us call this function $f(x,t^\star)$. Then $f(x,t)$ solves the heat equation in $I_{\delta}\times [t_\delta, t^\star]$.
Below we show that $F(x)$ is non-degenerate at $x=p(t)$, thus eliminating (ii). First let us show the analyticity of $F$ in $I_{\delta}$.

\begin{lemma}\label{lemma-analytic}
For every $t_\delta<t\leq t^\star$ fixed, the function $f(\cdot,t)$ is analytic in $I_\delta$. Moreover the function $f(\cdot,t^{\star})$ is positive to the left of $p(t^{\star})$, and negative to the right of  $p(t^{\star})$.\end{lemma}
\begin{proof}
We write
\begin{align*}
f(x,t) &= \int_{-1}^{+1} \Gamma(x,x^{\prime};t)f_I(x^{\prime})\;dx^{\prime} \\
&+
\int_0^t \left[\Gamma(x,(p(t^{\prime})-a);t-t^{\prime}) - \Gamma(x,(p(t^{\prime})+a);t-t^{\prime}) \right]\lambda(t^{\prime})
\;dt^{\prime}\end{align*}

The part of the solution coming from the linear term is analytic (see \cite{Evans} Section 2.3.c). For the nonlinear part, it is an easy computation to show that the function $f$ extended to the complex plane $f(x) \to f(z)$ with $z = x + iy$ is complex differentiable i.e., satisfies the equation 
$$
\frac{\partial f}{\partial x} + i \frac{\partial f}{\partial y}=0.
$$

As for $t=t^*$, note that the formula still holds for $F(x) = f(x,t^\star)$.

\medskip

For the second claim observe that $f(x,t)$ solves the heat equation in the domain $I_{\delta}\times [t_{\delta}, t^{\star}]$ with continuous zero set $p(t)$, and with initial data $f(x,t_{\delta})$ which was positive to the left of $p(t_{\delta})$ and negative to the right of $p(t_{\delta})$. Furthermore on the lateral boundary of $I_{\delta}$, $f(x,t)>0$ on the right and $f(x,t)<0$ on the left. Therefore the claim follows directly from the strong maximum principle for solutions of heat equation.
\end{proof}

We try now to take care of the possibility [ii]. Assume that
$$f_x(p(t^\star),t^\star) =0.$$
The following holds

\begin{lemma}\label{zero-derivatives}
Let $t^\star$ the first time at which the gradient at the free boundary degenerates, i.e $f_x(p(t^\star),t^\star) =0.$ Then, also
$$
f_{xx}(p(t^\star),t^\star) =0.
$$
\end{lemma}
\begin{proof}
The proof follows from a very simple observation:  
Suppose that $f_{xx}(p(t^\star),t^\star) >0$. This fact cannot coexists with the properties of the function $f$: $f$ is positive for $x<p(t^*)$ and negative for $x>p(^*)$ with purportedly vanishing derivative and very smooth across the free boundary (in fact $f$ is spatially analytical in space in a neighborhood of the free boundary, as we have shown).
\end{proof}

The question is now how to define the velocity of the free boundary at the time $t^\star$. Note in fact that as long as $f_x(p(t),t)$ stays away from zero, the speed of $p(t)$  is given by
$$\dot{p}(t) = - \frac{f_{xx}(p(t),t)}{f_x(p(t),t)},$$
which, as we have just shown, is an indefinite form if the first derivative degenerates at $x=p(t)$. 

Then, in order to find define the speed of $p(t)$ at this critical time $t^\star$, let us take the Taylor expansion of the equation $f(p +\Delta p , t+ \Delta t ) =0$. 
It holds
\begin{align*}
0 =& f(p +\Delta p , t+ \Delta t ) = f(p,t) + f_x(p,t)\Delta p + f_t(p,t)\Delta t \\
&+ \frac{1}{2!} [ f_{xx} \Delta p^2 + 2 f_{xt}\Delta p\Delta t + f_{tt} \Delta t^2 ] \\
& +  \frac{1}{3!} [ f_{xxx} \Delta p^3 + 2 f_{xxt}\Delta p^2\Delta t  +2 f_{xtt}\Delta p\Delta t^2  +  f_{ttt} \Delta t^2 ] + ...
\end{align*}
At the time $t^\star$, from the Taylor expansion we get that
$$
\dot {p} = -\frac{ f_{xxxx}}{2f_{xxx}}.
$$
We would like now do determine if this speed in finite or infinite: from the Hopf Lemma we know that in a parabolic problem, as long as the free boundary has finite speed, the derivative at the boundary is strictly positive (or negative) for all time. This implies that in our case, at the time $t^\star$ the speed of $p$ is not finite. Consequently the value $ f_{xxx}(p(t^\star),t^\star)$ has to be zero. 
Now we can repeat Lemma \ref{zero-derivatives} for the derivatives $f_{xxx}$ and $f_{xxxx}$ and prove that also
$$ f_{xxxx}(p(t^\star),t^\star)=0.$$ 
Again using the Taylor expansion we would then define
$$
\dot {p} = -\frac{ f_{xxxxxx}}{3f_{xxxxx}}, 
$$
etc.  Repeating Lemma \ref{zero-derivatives} iteratively, in is seen that all the derivatives of the function $f$ at the time $t^\star$ vanish.

Since all the derivatives of $f(x,t^\star)$ are equal to zero at $p(t^\star)$ from the previous lemma, analyticity implies that $f(x,t^\star)$ is zero in the entire neighborhood $I_{\delta}$ as well as the point $p(t^\star)$. \\
This is a contradiction since we have showed in  Lemma~\ref{lemma-analytic} that any solution at the time $t^\star$ is strictly positive for $x<p(t^\star)$ and negative for $x>p(t^\star)$.

Therefore the gradient of the function at the free boundary $p(t)$ stays strictly negative for any time $t\ge 0$. \\
This eliminates case [ii] and we are done if there are no puddles.  \\


\subsection{Puddle case: $\mathcal{L}:=\liminf_{t\to t^\star} p(t)< \mathcal{R}:=\limsup_{t\to t^\star} p(t)$}

Here the situation requires more careful analysis: first we claim that if a puddle forms at $t=t^\star$ then all spatial derivatives of $f$ are bounded in the whole domain, at least along a sequence of times converging to $t^\star$. (For precise statement see Proposition \ref{derivative} below). This result, along with an analyticity argument, will immediately yield a contradiction.  It is remarked that in the regularized case, such results follow, essentially from quadrature.  However, these estimates deteriorate as the regularization is removed so we shall not pursue this venue.  

In the case of a puddle, it is clear that
$$
\int_{\mathcal{L}}^{\mathcal{R}} |f|(x,t) dt \to 0 \hbox{ as } t\uparrow (t^\star).
$$
 (This is because the zero travels between the endpoints of the puddle so fast that there is not enough time for mass to diffuse into or flux through the zero set as $t\to t^\star$.)
 
We begin by choosing $t_0$ close enough to $t^\star$ such that 
$$
\int_{\mathcal{L}}^{\mathcal{R}} |f|(x,t) dt \leq \frac{1}{10}\min(M_b,M_p).\leqno \text{(BH)} 
$$

\subsubsection{Derivative bounds}

\begin{proposition}\label{derivative}
Suppose
$$
\liminf_{t\to t^{\star}} p(t) = \mathcal{L} < \mathcal{R}=\limsup_{t\to t^{\star}} p(t).
$$
Then there exists a constant  $C(m)$ such if $\beta > 0$ is sufficiently small then the following holds: for any sequence $t_n\to t^{\star}$ such that $p(t_n)\to \mathcal{R}$
 $$
\limsup_{n\to\infty} |f^{(m)}(x,t_n)| \leq C(m)\hbox{ for any } 
x\in [\mathcal R - \beta,\mathcal R + \beta].
 $$
A parallel result holds for $\mathcal{L}$.
\end{proposition}

\begin{proof}

Since $p(t_{n}) \to \mathcal R$, the number $\beta$ is chosen so that $\underline{a}(\mathcal R) > 3\beta$.  This means that for all $n$ sufficiently large, 
\begin{equation}
\label{CVG}
p(t_{n}) + \underline{a}(p(t_{n})) > \mathcal R + 2\beta
\end{equation}
i.e., during the times of specific interest the sink is well outside the interval under consideration.  Moreover at some point $\tilde t$ in 
$[\mathfrak t_{0}, t^{\star}]$, 
$p(t) - \underline{a}(p(t)) < \mathcal R - 2\beta$
for $t \geq \tilde t$.  Thus, for all times after $\tilde t$, for any 
$x_{0} \in [\mathcal R - \beta, \mathcal R + \beta]$, the source only represents a distant agitation.
With this in mind, let us reset $\mathfrak t_{0}$ such that all of the above holds for all $n$ and/or for all times greater than 
$\mathfrak t_{0}$.   Moreover, due to an up and coming plethora 
of indices, we might as well define $t_{n} =: t^{\#}$, with (only) the property of Eq.(\ref{CVG}) to be reserved for later.

It will prove convenient to work with with the auxiliary variable which measures the time remaining:
\begin{equation}
\theta:=\theta(t)  = t^{\#} - t
\end{equation}
for $t \in [\mathfrak t_{0}, t^{\#}]$.  Let $B(\theta)$ denote a slowly diverging function e.g., to be explicit
\begin{equation}
B = \sqrt{B_{0}}|\log\theta|^{1/2}
\end{equation}
will be adequate  with $B_{0}$ chosen for the purposes at hand.  For
$x_{0} \in [\mathcal R -\beta, \mathcal R + \beta]$ let us estimate magnitude of the $m^{\text{th}}$ derivative of $f(x_{0}, t^{\#})$.

We denote by $X_{0} :=x_0- \underline{a}$ which locates $p(t)$ when the source contributes to
$f(x_{0}, \cdot)$ and its derivatives.  If we denote by
$C(t)$ the distance between $p(t)$ to $X_{0}$ then by the Green's function formula
\begin{equation}
|f^{(m)}(x_{0}, t^{\#})|  \leq
\int_{\mathfrak t_{0}}^{t^{\#}}\frac{\lambda(t)}{\theta^{L}}\text{e}^{-C(t)^{2}/\theta}dt
\end{equation}
for some $L = L(m)$.  We start by defining the set
$\mathbb H \subset [\mathfrak t_{0}, t^{\#}]$
via
\begin{equation}
\mathbb H = \{ t \mid C(t) < B(\theta) \sqrt \theta \}.
\end{equation}
Our first claim is that $\mathbb H$ is the only important set for the
ostensible development  of singularities in $f(x_{0}, \cdot)$.  Indeed
\begin{equation}
\int_{\mathbb H^{c}}
\frac{\lambda(t)}{\theta^{L}}\text{e}^{-C(t)^{2}/\theta}dt
\leq
\int_{\mathbb H^{c}}
\frac{\lambda(t)}{\theta^{L}}\text{e}^{-B^{2}}dt
\leq
\int_{\mathfrak t_{0}}^{t^{\#}} \lambda(t) \theta^{B_{0} - L}dt
\end{equation}
which converges and is small independent of 
$t^{\#} < t^{\star}$ for $B_{0}$  large enough.
Henceforth we may focus on events that take place when $t \in \mathbb H$ where, for all intents and purposes, there is no help from the exponential factors.  (But, on the positive side, $\mathbb H$ is disjoint from the tail end of $[\mathfrak t_{0}, t^{\#}]$)

To aid with our objectives, it will be convenient to divide
$[\mathfrak t_{0}, t^{\#}]$ into disjoint regions that are of equal size on a logarithmic scale:  Let $H \gg 1$ denote a sufficiently large number the precise (minimum) value of which will be determined in what is to follow.
Roughly speaking, we wish the $k^{\text{th}}$ region to be of size $H^{-1}$ of the $k-1^{\text{st}}$.  Specifically, we may proceed as follows:  The $k^{\text{th}}$ region will be denoted by $g_{k}$, $k = 0, 1, \dots$ and the size, $|g_{k}|$, will satisfy
$|g_{k}| = H^{-k}|g_{0}|$.
(Thus $|g_{0}| = \frac{H-1}{H}(t^{\#} - \mathfrak t_{0})$). Therefore
\begin{equation}
g_{0}  =  [\mathfrak t_{0}, \mathfrak t_{0} + |g_{0}|),
\hspace {4 pt} \dots \hspace{4 pt},
g_{k}  =  [\mathfrak t_{0}  + |g_{0}| + \dots + |g_{k-1}|,
 \mathfrak t_{0} + |g_{0}| + \dots + |g_{k}|].
\end{equation}
Note that the value of $\theta$ at the right end of $g_{k}$ is a constant (that is very near one) times  the size of $g_{k+1}$.

We shall also need two spatial regions:  $\mathbb A_{2}$, $\mathbb A_{3}$ which are given by
$\mathbb A_{2} =  (X_{0} + b, X_{0} + 2b)$
$\mathbb A_{3} =  (X_{0} + 2b, X_{0} + 3b)$
with $b$ chosen small enough so that these sets are inside
$[\mathcal L, \mathcal R]$ lying well (on the scale of $b$) to the left of $\mathcal R$.  Moreover $b$ is large enough 
(or $\mathfrak t_{0}$
late enough) so that $b$ is much larger than 
$B(\theta(\mathfrak t_{0}))\sqrt{\theta(\mathfrak t_{0})}$.
We shall define {\it epochs} that are punctuated by certain
exits from and
entrances to the region
$\mathbb A_{3}$.
The beginning of an epoch, denoted by $\overline{\tau}_{\text{min}}$ is when $p(t)$ enters $\mathbb A_{2}$ from $\mathbb A_{3}$ and will not revisit
$\mathbb A_{3}$ before first having touched the appropriate
$B\sqrt\theta$ neighborhood of $X_{0}$.  Let $\overline{\tau}_{\text{min}}$ be the first time during this epoch where $p(t)$ visits $X_0+ B\sqrt{\theta}$.
The time $\tau_{\text{max}}$
is when $p(t)$ leaves the ostensibly smaller $B\sqrt\theta$ neighborhood
and does not touch these neighborhood types till another visit to 
$\mathbb A_{3}$. To be specific, $\tau_{\text{max}}$ will be the moment of this departure, so 
$p(\tau_{\text{max}})=X_0+ [B\sqrt{\theta}]|_{\tau_{\text{max}}}$.
There is a third time, namely $\tau_{\text{min}}$ when $p(t)$ actually enters the
appropriate $B\sqrt\theta$ neighborhood but this time does not play a particularly major
r\^{o}le.

Our first claim is that if epochs are localized e.g., to a single $g_{k}$ then their contribution
to $||f^{(m)}(x_{0},\cdot)||$ is tractable.  Indeed, assume for simplicity that times
$\overline{\tau}_{\text{min}}^{[1]} < \tau_{\text{min}}^{[1]} <\tau^{[1]}_{\text{max}} < \dots <
\overline{\tau}_{\text{min}}^{[J]} < \tau_{\text{min}}^{[J]} < \tau^{[J]}_{\text{max}}$
are the only punctuation marks in $g_{k}$.
Then
\begin{equation}
\label{HHK}
\int_{\tau_{\text{min}}^{[1]}}^{\tau_{\text{max}}^{[J]}}
\hspace{0 pt}
\frac{\lambda(t)}{\theta^{L}}\text{e}^{-C(t)/\theta^{2}}\chi_{\mathbb H}dt
\leq
\frac{1}{\theta_{k+1}^{L}}\int_{g_{k}\cap \mathbb H}
\hspace{-10 pt}
\lambda(t)dt
\end{equation}
where $\theta_{k+1} := t^\#-(t_0+...+|g_k|) \approx |g_{k+1}|$ is the time remaining by the end of
$g_{k}$.  

Next we show that the integral to be done is actually of the order of
$\exp[-1/\theta_{k}]$. This will be facilitated by the perspective of $\rho_{p}$:  Note that at
$t = \overline{\tau}_\text{min}^{[1]}$,
the density $\rho_{p}(x,\overline{\tau}_\text{min}^{[1]})$ vanishes
for $x < X_{0} + 2b$ and has some positive profile for
$x > X_{0} + 2b$.  Starting from this profile, we are supposed to compute the flux of $\rho_p$ through
$p(t)$ while $t\in \mathbb H$. It is not hard to see
that this is less than the total flux through (the vicinity of)
$X_{0}$ with Dirichlet boundary conditions at $x=X_0+B\sqrt\theta$, and an initial profile of a
delta mass at $X_{0} + 2b$.  
The strength of the $\delta $--mass might initially taken to be 
the total of $M_{p}$.  But to account for the possibility of reflux from the $p$--source (sink) -- which is much further away, let us add to this the quantity 
$$
\Lambda := \int_{0}^{t^{\star}}\lambda dt < \infty
$$ 
which is all that ever has fluxed and ever will flux.

This leads to the estimate
\begin{equation}
\int_{g_{k}\cap \mathbb H}\lambda(t)dt
\hspace{2 pt}\leq \hspace{2 pt}
C_{2}\text{e}^{-c_2/\theta_{k}}
\end{equation}
where $C_{2} < \infty$ and $c_{2} > 0$ are constants which do not depend on $k$, or the end time  $t^{\#}$.

The key issue is therefore to show that the epochs do not extend over many $g$--scales.
Indeed, it is remarked, an extension of the above reasoning tracking a single epoch through scales $\kappa \leq k \leq K$ would yield a bound of the form
\begin{equation}
\label{KK}
\frac{\tilde{C}_{1}}{\theta_{K}^{L}}\text{e}^{-c_{1}/\theta_{\kappa}}
\end{equation}
replacing the corresponding right hand side of Eq.(\ref{HHK}).
And, as unlikely as it may seem, if $K \gg \kappa$ this could be large.
We further remark that some positive powers of $C(t) = |X_{0} - p(t)|$
originating from various places
are available (inside the integrand) for the estimate.  But even in the best case scenario -- namely the estimate for the norm of $f(x_{0},t)$
itself -- there is still a logarithmic divergence ($\propto K - \kappa$) weighing in against the exponential prefactor.  
Finally, it might be {\it technically} noted that since 
$t^{\#} < t^{\star}$ it must be the case that all but a finite number of the $g_k$ are devoid of epochs.  However the technicality  
is of no practical significance since we seek bounds that are uniform in $t^{\#} < t^{\star}$.
We turn to the task at hand.

\medskip

Consider an epoch defined by the times
$\bar{\tau}_{\text{min}}$ and $\tau_{\text{max}}$.  Let $\mathbb A_{3,2}$ denote the left and middle third of the
$\mathbb A_{3}$ region namely
\begin{equation}
\mathbb A_{3,2}  =
\{x\mid X_{0} + 2b \leq x \leq X_{0} + \frac{8}{3}b \}
\end{equation}
and let us estimate the accumulation of mass
in the region $\mathbb A_{3,2}$ at time $t = \tau_{\text{max}}$.  
Note that in addition to the flux from the outside of the region, 
there is the initial mass that was in $\mathbb A_{3}$
at time $t = \bar{\tau}_{\text{min}}$.  This is actually of 
{\it negative} utility.  We cannot rely on it staying 
in the region and we do not wish to account for where it might go 
during $\bar{\tau}_{\text{min}} < t < \tau_{\text{max}}$.  
However, in accord with condition BH, this could only account for
$10\%$ of $M_{p}$. Then by mass conservation there is  a significant portion of $M_p$ outside of 
to the right of $\mathcal{R}$.

Thus, we shall rely on this material that diffusing in from the right.
Here, it is worthwhile to recollect that the underlying condition of the epoch is that
$p(t)$ stay out of $\mathbb A_{3}$.  So, placing the guaranteed fraction of $M_{p}$ at the extreme right -- $x = +1$ -- and placing Dirichlet boundary conditions at $x= X_{0}+2b$ we obtain the estimate
\begin{equation}
\label{middle}
Q_{3,2}  =  \int_{\mathbb A_{3,2}}\rho_{p}(x,\tau_{\text{max}})dx
\geq  C_{3,2}\text{e}^{-c_{3,2}/[\theta(\overline{\tau}_{\text{min}}) - \theta(\tau_{\text{Max}})]}
\end{equation}
where $C_{3,2} > 0$ and $c_{3,2} < \infty$ are constants independent of the parameters of the epoch.

Now just about all of $Q_{3,2}$ and more will be swept into $p(t)$ between time $\tau_{\text{max}}$ and the next time that $p(t)$ crosses all the way to the right side of the region $\mathbb A_{3}$.  This may happen immediately but it might incur the passage of other possible (long) epochs.  However, from the stipulations about the time sequence $(t_{n})$ it is 
inevitable that this will happen before
$t = t^{\#}$.

Our final substantive claim is that
$Q_{3,2}$ is essentially dominated by the flux through the left boundary of $\mathbb A_{2}$ in the time
interval $[\tau_{\text{max}}, t^{\#}]$ by a delta--source of strength $M_{p}+M_{b} + \Lambda$
placed at $X_{0} + B\sqrt{\theta}|_{\tau_{\text{max}}}$ at time $\tau_{\text{max}}$.  The implication provides an upper bound which we state as a separate lemma:

\begin{lemma}\label{upperbound}
Let $Q_{3,2}$ be as described. Then
\begin{equation}
Q_{3,2}  \leq  C_{4}\text{e}^{-c_{4}/\theta(\tau_{\text{Max}})}
\end{equation}
where $c_4$ and $C_4$ do not depend on the parameters of the epoch nor on the time $t^{\#}$.
\end{lemma}

On the basis of Lemma \ref{upperbound} and Eq.(\ref{middle}), we may conclude that
\begin{equation}
\label{DZ}
\theta(\tau_{\text{max}}) \geq c_{1}\theta(\tau_{\text{min}})
\end{equation}
with $c_{1} > 0$ independent of the parameters of the epoch
and the time $t^{\#}$.  
Thus consider $H$ to be sufficiently large ($H>1/c_1$) and let us double cover $[\mathfrak t_{0}, t^{\#}]$
with two overlapping partitions of the general type described -- both with parameter $H$ -- in such a way that (on the logarithmic scale) the midpoints of one partition form the endpoints of the other.

Then according to Eq.(\ref{DZ}), all epochs are caught in individual elements of one or the other partition.  This enables us to sum (twice) the analogue of Eq.(\ref{KK}) with matching $k$'s all the way to $k = \infty$
which demonstrates the desired result.
\end{proof}

\noindent \textbf{Proof of Lemma~\ref{upperbound}}.

\medskip

Consider the situation at $t = \tau_{\text{max}}$ where that there is non--zero $\rho_{b}$ for $x < X_{0}$ and non--zero $\rho_{p}$ for $x > X_{0}$.
The initial densities (and ensuing fluxes) are, evidently such that soon enough, $p(t)$ will exit from the neighborhood of $X_{0}$ and cross all the way into $\mathbb A_{3}$.  Later -- maybe much later on the log--scale -- it will cross all the way {\it through} $\mathbb A_{3}$

Now consider $\tilde{\rho}(x,t)$ defined in the domain $\Sigma^{\#}:=[-1, X_0+3b]\times [\tau_{\text{max}}, t^\#]$ as follows:
\begin{itemize}
\item[(a)]  $(\partial_t-\Delta)\tilde{\rho}= 
\lambda(t)\delta_{x=p(t) - \underline{a}}$ ; \\
\item[(b)] $\partial_x\tilde{\rho}(-1,t)=0$ ;\\
\item[(c)] $\tilde{\rho}(X_0+3b,t) \equiv 0$ ;\\
\item[(d)] $\tilde{\rho}(x,\tau_{\text{max}})=\rho_b(\cdot,\tau_{\text{max}})+\varphi(x)-\rho_p(x,t)\chi_{(X_0+2b, X_0+\frac{8}{3}b)}(x)$ ;\\
\end{itemize}
\noindent where $\varphi$ is the mirror image of $\rho_p(\cdot,\tau_{\text{max}})$ in $\mathbb A_{3}$ reflected about $x=X_0+b$, i.e.,
$$
\varphi(x) := \rho_p(2(X_0+b)-x,\tau_{\text{max}})\chi_{(X_0-b, X_0)}(x).
$$

  We denote $\alpha(t):=\{x: \alpha(x,t)=0\}$ and denote by $t^{\ddag}$ be the first time when $\alpha(t)$ touches $X_0+3b$.
  There are two observations:  
  
  \noindent $\bullet$ First, by maximum principle for caloric functions it follows that $\rho\leq \tilde{\rho}$ in $\Sigma^{\#}$, and thus
  \begin{equation}
  p(t) \leq \alpha(t) \quad \hbox{ for } \tau_{\text{max}} \leq t\leq t^{\ddag}
  \end{equation}
 This, of course, places $t^{\ddag} \leq t^{\#}$.

\noindent $\bullet$ Secondly, observe that $\tilde{\rho}(x,t) \geq \tilde{\rho}(2(X_0+b)-x,t)$ in $\Sigma$, again by the maximum principle for caloric functions. In particular,
\begin{equation}
\alpha(t)  \in [X_0+b,X_0+3b] = \mathbb A_{2}\cup \mathbb A_{3};
\hspace{8 pt}
 \tau_{\text{max}} \leq t \leq t^{\ddag}.
\end{equation}
Since $\int \tilde{\rho}(x,\tau_{\text{max}})dx=Q_{2,3}$ and $\int \tilde{\rho}(x,t^{\ddag})dx=0$,
$t = t^{\ddag}$, it follows that
\begin{equation}
\label{flux}
Q_{2,3}= \int_{\tau_{\text{max}}}^{t^{\ddag}}(-\tilde{\rho}_x)(\alpha(t'),t') dt' + \int_{\tau_{\text{max}}}^{t^{\ddag}}\tilde{\rho}_x(X_0+3b,t') dt'.
\end{equation}

Moreover since we know
$\alpha_{t} \in \mathbb A_{3}\cup \mathbb A_{2}$ for all $t \in [\tau_{\text{max}}, t^{\ddag}]$,
by obvious dominance the second term in Eq. (\ref{flux}) is bounded by the flux of $h(x,t)$ through $x=X_0+3b$, where $h(x,t)$ solves the heat equation in the region to the right of $X_{0} + b$ during the times $[\tau_{\text{max}},t^{\ddag}]$
with initial data the same as $\tilde{\rho}$ at $t=\tau_{\text{max}}$ and Dirichlet conditions both sides of 
$\mathbb A_{2}\cup \mathbb A_{3}$.
Since  $t^{\ddag} < t^{\#}$ we may write that the {\it fraction} of $Q_{3,2}$ that is lost via the right boundary is less than
$C_{5}\text{e}^{-c_{5}/\theta(\tau_{\text{max}})}$ where the $C_5$ and $c_5$ only depends on $b$.  (So that the {\it total} which is lost is no more than $Q_{2,3}C_{5}\text{e}^{-c_{5}/\theta(\tau_{\text{max}})}$.
As far as we are concerned it is sufficient that this is less than half of $Q_{2,3}$.)

It remains to estimate the first term in Eq.(\ref{flux}).  For this purpose, let us consider $g(x,t)$ solving the same problem $(a)-(d)$ in $\Sigma$ as $\tilde{\rho}$, 
except that $g(x,t)$ has Dirichlet boundary conditions at $x = X_{0} + b$.  By caloric inequalities the flux of 
$\tilde \rho$
through 
$\alpha(t)$ (where $\alpha(t) \geq b$) in the interval 
$[\tau_{\text{max}},t^{\ddag}]$
is, from the perspective of the left, less than the flux of 
$g$ through the line $x = X_{0} + b$.  
Finally we may further modify $g$ to a $\tilde g$
which has the Dirichlet condition at $x = X_{0}+ b$ 
has all conceivably available mass 
-- namely $M_{p} + M_{b}$ -- is placed at 
$X_{0}$.  Moreover, the source term is placed as far forward as possible and allowed (more than) all of its available flux as soon as possible.  This amounts to adding $\Lambda$ to the $M_{p}\ \& \ M_{b}$ which are in the mass at $x = X_{0}$
Thence

$$
\begin{array}{lll}
 \int_{\tau_{max}}^{t^{\ddag}}(-\tilde{\rho}_x)(\alpha(t'),t') dt'&\leq& \int_{\tau_{max}}^{t^{\ddag}} (-\tilde{g}_x)(2b,t') dt'\\ \\
&\leq & C_{4}e^{-c_{4}/\theta(\tau_{max})}
\end{array}
$$
where $c_{4}$ and $C_{4}$ only depends on $b$, $M_{p}$, etc.
\hfill$\Box$

\medskip

Now may now we finish the proof of Theorem~\ref{thm:main}.

\begin{corollary}
There does not exist  a finite blow-up time $t^\star$.
\end{corollary}

\begin{proof}
We already eliminated the non-puddle case in Section 3.1, so let us discuss the puddle case.
Let $\mathcal{L}$ and $\mathcal{R}$ be as above. Let  $t_n\to t^\star$ be such that  $p(t_n)\to \mathcal{R}$.
In particular, choose $n$ large enough such that $|p(t_n) - \beta | < \epsilon \ll \frac{1}{8}\underline{a}$.
Call $f_n(x) = f(x,t_n)$. The function $f_n(x)$ has a unique zero at $p(t_n)$. In Lemma \ref{lemma-analytic} it was proven that (for an sufficiently large) $f_n$ is an analytic function in any subset of the open interval $(p(t_n)-\frac{1}{8}\underline{a}, p(t_n) +\frac{1}{8}\underline{a})$. This implies that $f_n(x)$ is analytic in the interval 
$(\mathcal R - \beta, \mathcal R + \beta)$ for large enough $n$.

Moreover due to Proposition~\ref{derivative} the sequence $\{f_n\}$ is uniformly bounded and uniformly Lipschitz continuous in the $2\beta$ 
neighborhood of $\mathcal{R}$. Hence Ascoli-Arzela Theorem 
ensures the existence of a subsequence $f_{n_k}$ such that $f_{n_k} \to \phi(x)$  as $t_n\to t^\star$ uniformly in  $(\mathcal{R}-\beta, \mathcal{R}+\beta)$.  \\


Now, since the sequence of analytic functions $f_{n_k}$ converges uniformly, the limiting function $\phi$ is also analytic in $(\mathcal{R}-\beta, \mathcal{R}+\beta)$. Now let us observe the profile of $\phi$ near $\mathcal{R}$. From the choice of the sequence $t_n$, it is clear that $\phi(x)$ is positive to the right  of $x=\mathcal{R}$. 
(Indeed for any $\e>0$,   for sufficiently large $n$ depending on $\e$ the function $f(x,t)$ solves the heat equation with source term in $\Sigma_\e:=[-\mathcal R + \epsilon, +1]\times [t_n,t^\star)$ with positive boundary data at $\mathcal R + \epsilon$ and $t=t_n$ and Neumann boundary data at $x = 1$. Therefore $f(x,t)$ stays strictly positive in $\Sigma_\e$, staying uniformly away from zero as $t\to t^\star$. )
On the other hand $\phi$ is identically zero to the left of 
$\mathcal R$. These two facts interdict the possibility of analyticity for the function $\phi$.
\end{proof}

\noindent \textbf{Acknowledgements: }
L. Chayes is supported by the NSF under the grant DMS--0805486.
M.d.M. Gonzalez is supported by Spain Government project MTM2008-06349-C03-01, and NSF DMS-0635607.
M.P. Gualdani is supported by NSF-DMS 0807636 as well as the grant DMS-0306167.  I. Kim is supported by NSF DMS-0700732 and a Sloan Foundation fellowship.


\bibliographystyle{abbrv}

\end{document}